\documentclass[12pt,leqno,openany]{article}
\usepackage{amsmath,amssymb,amsthm,textcomp,times}
\usepackage[T1]{fontenc}
\usepackage[latin1]{inputenc}
\usepackage{colortbl}
\usepackage[english]{babel}
\usepackage{dsfont}
\newtheorem{theorem}[subsection]{Theorem}

\newtheorem{definition}[subsection]{Definition}
\newtheorem{lemma}[subsection]{Lemma}
\newtheorem{proposition}[subsection]{Proposition}

\newtheorem{remark}[subsection]{Remark}

\def \R{ \mathbb R }

\def \E{ \mathbb  E }

\def \P{ \mathbb P  }

\begin{document}
\title{KPZ formula for log-infinitely divisible multifractal random measures}
\maketitle
\begin{center}
{R\'emi Rhodes, Vincent Vargas \\
\footnotesize

 CNRS, UMR 7534, F-75016 Paris, France \\
  Universit{\'e} Paris-Dauphine, Ceremade, F-75016 Paris, France} \\

{\footnotesize \noindent e-mail: \texttt{rhodes@ceremade.dauphine.fr},\\
\texttt{vargas@ceremade.dauphine.fr}}
\end{center}

\begin{abstract}
We consider the continuous model of log-infinitely divisible multifractal random measures (MRM) introduced in \cite{bacry}. If $M$ is a non degenerate multifractal measure with associated metric $\rho(x,y)=M([x,y])$ and structure function $\zeta$, we show that we have the following relation between the (Euclidian) Hausdorff dimension ${\rm dim}_H$ of a measurable set $K$
and the Hausdorff dimension ${\rm dim}_H^{\rho}$ with respect to $\rho$ of the same set: $\zeta({\rm dim}_H^{\rho}(K))={\rm dim}_H(K)$. Our results can be extended to higher dimensions in the log normal case: inspired by quantum gravity in dimension $2$, we consider the $2$ dimensional case.

\end{abstract}
\vspace{1cm}
\footnotesize


\noindent{\bf Key words or phrases:} Random measures, Hausdorff dimensions, Multifractal processes.


\normalsize


\section{Introduction}
Multiplicative cascades are random measures that were introduced by Mandelbrot in \cite{cf:Man} to model the energy dissipation of a turbulent flow. This model, which arises as the limit of discrete random multipliers, has been the object of numerous studies in probability theory (see for instance \cite{cf:Liu} for an account on the achieved results). In the beautiful note \cite{Benj}, inspired by the work of \cite{cf:DuSh}, the authors related the Hausdorff dimension ${\rm dim}_H$ of a measurable set $K$  to the Hausdorff dimension of the same set in the random metric induced by the multiplicative cascade: this gave the so called KPZ formula in analogy with a similar formula in quantum gravity (\cite{cf: KPZ}).

In this work, we derive a similar formula in the context of log-infinitely divisible multifractal random measures (MRM) introduced by the authors in \cite{bacry}. MRM are scale invariant generalisations of the log normal model introduced in \cite{cf:Manturbulence} (and rigorously defined mathematically by Kahane in \cite{cf:Kah}) and the log Poisson model studied in \cite{Bar}. MRM have been used as models of the energy dissipation in a turbulent flow (see \cite{cf:Fr}) and of the volatility of a financial asset (see \cite{cf:BaKoMu}, \cite{cf:DuRoVa}); as such, MRM  are much more realistic models than multiplicative cascades whose construction relies on a discrete dyadic decomposition of the unit interval. In particular, this dyadic dependent construction entails that multiplicative cascades have non stationary increments which is not the case of MRM.

The following note is organized as follows: section 2 reminds the definition and main properties of MRM.
 Section 3 reminds the background on Hausdorff dimensions needed in the proof of the main theorem. In section 4, we state the main theorem in dimension 1: theorem \ref{main}. In section 5, we give the 2-dimensional analog for MRM and the Gaussian free field (inspired by quantum gravity). In section 6, we give the detailed proof of theorem \ref{main}: our proof follows tightly the one given in \cite{Benj} for multiplicative cascades. Nevertheless, the main estimates needed to carry out the proof are more difficult for MRM (the use of scale invariance is crucial: see item 4. in proposition 2.5 below). In section 7, we prove the theorems of section 5.

\begin{remark}
At the time we write this article, we have not seen the work of Duplantier and Sheffield (\cite{cf:DuSh}) which inspired the note \cite{Benj}: we are therefore indirectly indebted to them. It seems that in \cite{cf:DuSh} the authors prove a result similar to our theorem \ref{th:GFF} (see below) using the theory of large deviations for Gaussian processes: it would be interesting to compare their result with our theorem \ref{th:GFF}. In this article, we do not use large deviation theory; we prove theorem \ref{th:GFF} by a straightforward adaptation of the proof of theorem \ref{main} (valid in dimension 1 for log infinitely divisible measures and in particular for log Gaussian measures).  
\end{remark}



\section{Introductory background about MRM}

The reader is referred to \cite{bacry} for all the proofs of the results stated in this section.

\vspace{2mm}
{\bf Independently scattered infinitely divisible random measure.} Let $S^+$ be the half-plane
$$ S^+=\{(t,y);t\in\R,y\in\R_+^*\}$$ with which we associate the measure (on the Borel $\sigma$-algebra $\mathcal{B}(S^+)$)
$$\theta(dt,dy)=y^{-2}dt\,dy. $$
The characteristic function of an infinitely divisible random variable $X$ can be written as $\E[e^{iqX}]=e^{\varphi(q)}$, where $\varphi$ is characterized by the L\'evy-Khintchine formula
$$ \varphi(q)=i mq-\frac{1}{2}\sigma^2q^2+\int_{\R^*} (e^{iqx} -1 -iq\sin(x))\,\nu(dx)$$
and $\nu(dx)$ is the so-called L\'evy measure. It satisfies $\int_{\R^*}\min(1,x^2)\,\nu(dx)<+\infty$.

Following \cite{bacry}, we consider an independently scattered infinitely divisible random measure $\mu$ associated to $(\varphi,\theta)$ and distributed on the half-plane $S^+$ (see \cite{Ros}).
More precisely, $\mu$ satisfies:


1) For every sequence of disjoint sets $(A_n)_n$ in $\mathcal{B}(S^+)$, the random variables $(\mu(A_n))_n$ are independent and
$$\mu\big(\bigcup_nA_n\big)=\sum_n\mu(A_n) \text{ a.s.},$$

2) for any measurable set $A$ in $\mathcal{B}(S^+)$, $\mu(A)$ is an infinitely divisible random variable whose characteristic function is
$$\E(e^{iq\mu(A)})=e^{\varphi(q) \theta(A)}. $$

We stress the fact that $\mu$ is not necessarily a random signed measure. Let us additionnally mention that there exists a convex function $\psi$ defined on $\R$ such that for all non empty subset $A$ of $S^+$:

- $\psi(q)=+\infty$, if $\E (e^{q\mu(A)})=+\infty$,

-$\E(e^{q \mu(A)})=e^{\psi(q)\theta(A)}$ otherwise. \\
Let $q_c$ be defined as $q_c=\sup \{q\geq 0; \psi(q)<+\infty\}$. For any $q \in[0,q_c[$, $\psi(q) < +\infty$ and $\psi(q)=\varphi(-iq)$.

\vspace{2mm}
{\bf Multifractal Random Measures (MRM).} We consider an independently scattered infinitely divisible random measure $\mu$ associated to $(\varphi,\theta)$ such that $q_c>1$, namely that:
$$\exists\epsilon>0, \,\,\psi(1+\epsilon)<+\infty,$$ and $\psi(1)=0$.

\begin{definition}{\bf Filtration $\mathbf{\mathcal{F}_l} $.}
Let $\Omega$ be the probability space on which $\mu$ is defined. $\mathcal{F}_l$ is defined as the $\sigma$-algebra generated by $\{\mu(A); A \subset S^+, {\rm dist}(A,\R^2\setminus S^+)\geq l\}$.
\end{definition}

Let us now define the function $f:\R_+\to \R$ by
$$ f(l)=\left\{
\begin{array}{ll}
 l, &  \text{ if } l\leq T   \\
 T &  \text{ if } l\geq T
\end{array}\right.$$ The cone-like subset $A_l(t)$ of $S^+$ is defined by
$$ A_l(t)=\{(s,y)\in S^+; y\geq l, -f(y)/2\leq s-t\leq f(y)/2\}.$$
For forthcoming computations, we stress that $\theta(A_l(t))=\int_l^{+\infty}f(y)y^{-2}\,dy<+\infty$ and, for $l\leq T$, $ \theta(A_l(t))=\ln(T/l)+1$.

\begin{definition}{\bf $\mathbf{\omega_l(t)} $ process.}
The process $\omega_l(t)$ is defined as $\omega_l(t)=\mu(A_l(t))$.
\end{definition}

\begin{definition}{\bf $\mathbf{M_l(t)} $ measure.}
For any $l>0$, we define the measure $M_l(dt)=e^{\omega_l(t)}\,dt$, that is
$$ M_l(I)=\int_Ie^{\omega_l(r)}\,dr$$ for any Lebesgue measurable subset $I\subset \R$.
\end{definition}

\begin{definition}\label{wcmrm}{\bf Multifractal Random Measure (MRM).} With probability one,  there exists a limit measure (in the sense of weak convergence of measures)
$$M(dt)=\lim_{l\to 0^+}M_l(dt).$$
This limit is called the Multifractal Random Measure. The scaling exponent of $M$ is defined by
$$\forall q \geq 0, \quad \zeta(q)= q-\psi(q). $$
\end{definition}

\begin{proposition}{\bf Main properties of the MRM.}\label{propMRM}
\begin{enumerate}
\item the measure $M$ has no atoms in the sense that $M(\{t\})=0$ for any $t\in\R$.
\item The measure $M$ is different from $0$ if and only if there exists $\epsilon>0$ such that $ \zeta(1+\epsilon)>1$; in that case, $\E(M([0,t]))=t$.
\item if $\zeta(q)>1$ then $\E[M([0,t])^q]<+\infty $.
\item For any fixed $\lambda\in]0,1]$ and $l\leq T$, the two processes $(\omega_{\lambda l}(\lambda t))_{0\leq t \leq T}$ and $(\Omega_\lambda+\omega_l(t))_{0\leq t \leq T}$ have the same law, where $\Omega_\lambda$ is an infinitely divisible random variable independent from the process $(\omega_l(t))_{0\leq t \leq T} $ and its law is characterized by $\E[e^{iq \Omega_\lambda}]=\lambda^{-\varphi(q)}$.
\item For any $\lambda\in ]0,1]$, the law of the process $ (M([0,\lambda t]))_{0\leq t\leq T}$ is equal to the law of $(W_\lambda  M([0, t]))_{0\leq t\leq T}$, where $W_\lambda=\lambda e^{\Omega_\lambda}$ and $\Omega_\lambda$ is an infinitely divisible random variable (independent of $ (M([0, t]))_{0\leq t\leq T}$) and its characteristic function is
$$ \E[e^{iq \Omega_\lambda}]=\lambda^{-\varphi(q)}.$$
\item If $\zeta(q)\not =-\infty $ then
$$ \E\big[M([0,t])^q\big]=(t/T)^{\zeta(q)}\E\big[M([0,T])^q\big].$$
\end{enumerate}
\end{proposition}

\begin{proposition}{\bf Main properties of the scaling exponent.}
If there is $\epsilon>0$ such that $\zeta(1+\epsilon)>1$, the function $q\in[0,1]\mapsto \zeta(q)$ is continuous, strictly monotone increasing and maps $[0,1]$ onto $[0,1]$.
\end{proposition}

\section{Hausdorff dimension}
In this section, we just set out the minimal required background about the Hausdorff dimension to understand our main result and its proof. We refer to \cite{Falc} for an account on Hausdorff dimensions.
\begin{definition}
Let $(X,d)$ be a metric space. If $K\subset X$ and $s \in [0,+\infty[$, the $s$-dimensional {\it Hausdorff content} of $K$ is defined by
$$C^s_H(K)=\inf\left\{\sum_ir_i^s;\text{there is a cover of K by balls with radii } r_i>0\right\}.$$
Using the standard convention $\inf \emptyset=+\infty$, the {\it Hausdorff dimension} of $K$ is defined by
$${\rm dim}_H(K)=\inf \left\{s \geq 0; C^s_H(K)=0\right\}.$$
\end{definition}

\begin{lemma}{\bf (Frostman)}
Let $(X,d)$ be a metric space.The $s$-capacity of a Borelian set $K\subset X$
$$ {\rm Cap}_s(K)=\inf\left\{\Big(\int_{K\times K}|y-x|^{-s}\gamma(dx)\gamma(dy)\big)^{-1};\,\gamma \text{ is a Borel measure such that }\gamma(K)=1\right\}$$ is linked to the Hausdorff dimension of $K$ by the relation
$${\rm dim}_H(K)=\sup \left\{s\geq 0;  {\rm Cap}_s(K)>0\right\}.$$
\end{lemma}
\section{KPZ formula in one dimension}
If we define for $x,y\in\R$, $\rho(x,y)=M([x,y])$, then $\P$ a.s. $\rho$ is a random metric on $\R$. The interval $[0,T]$ can be seen as a metric space when it is equipped either with the Euclidean metric $|\cdot|$ or with the random metric $\rho$. The main purpose of this paper is to establish a relation between the Hausdorff dimension of a measurable set $K\subset [0,T]$  equipped with the Euclidean metric and its Hausdorff dimension with respect to the (random) metric space $([0,T],\rho)$.
\begin{theorem}\label{main}
Assume there is $\epsilon>0$ such that $\zeta(1+\epsilon)>1$ and that for all $q \in [0,1]$ we have $\psi(-q)< \infty$. Let $K\subset [0,T]$ be some deterministic and measurable nonempty set and $\delta_0$ its Hausdorff dimension with respect to the Euclidian metric. Then the Hausdorff dimension ${\rm dim}_H^\rho(K)$ of $K$ with respect to the random metric $\rho$ coincides $\P$ a.s. with the unique solution $\delta$ in $[0,1] $ of the equation $\delta_0=\zeta(\delta)$.
\end{theorem}

\begin{remark}
We can see $\rho$ as a strictly increasing function on $[0,T]$: $x \rightarrow \rho(0,x)$.
By definition of ${\rm dim}_H^\rho$ , we have $\P$ a.s.:
\begin{equation*}
\forall K \in \mathcal{B}(\rho([0,T])), \quad {\rm dim}_H^\rho(\rho^{-1}(K))={\rm dim}_H(K)
\end{equation*}
Applying the above equality to $\rho(K)$, we get an equivalent formulation to theorem \ref{main}: if $K$ is some deterministic measurable set, we get $\P$ a.s.:
\begin{equation*}
\zeta({\rm dim}_H(\rho(K)))={\rm dim}_H(K)
\end{equation*}
\end{remark}

\section{KPZ formula in 2 dimensions}
In this section, inspired by the KPZ formula in continuum quantum gravity (\cite{cf: KPZ}),  we consider the natural extension in dimension $2$ of the results of the previous section in the log normal case (the results of section 5.1 have analogs in all dimensions).

\subsection{The log normal MRM measure in dimension 2}
The log normal MRM in dimension 2 is the random measure $M$ in $\R^2$ defined formally by:
\begin{equation*}
\forall A \in \mathcal{B}(\R^2), \quad M(A)=\int_{A}e^{X(x)-\frac{1}{2}\E[X(x)^2]}dx
 \end{equation*}
where $(X(x))_{x \in \R^2}$ is a "Gaussian field" whose covariance is given by:
\begin{equation*}
\E[X(x)X(y)]=\gamma^2\ln^+\frac{R}{|x-y|}.
\end{equation*}
where $\gamma^2$ and $R$ are two positive parameters.
To give a rigorous meaning to $M$, one can use the theory of Gaussian multiplicative chaos introduced by Kahane in \cite{cf:Kah} or it's extension defined in \cite{cf:RoVa}. In this framework, the measure $M$ is the multiplicative chaos associated to the function $\ln^+\frac{R}{|x|}$ and it can be defined almost surely (see example 2.3 in \cite{cf:RoVa}) as the limit (in the space of Radon measures) as $l$ goes to $0$ of the random measures $M_{l}(dx)$ defined by:
\begin{equation*}
\forall A \in \mathcal{B}(\R^2), \quad M_{l}(A)=\int_{A}e^{X_{l}(x)-\frac{1}{2}\E[X_{l}(x)^2]}dx
 \end{equation*}
where $(X_{l}(x))_{x \in \R^2}$ is as centered Gaussian field whose covariance is given by:
\begin{equation*}
\E[X_{l}(x)X_{l}(y)]
=
\begin{cases}
\gamma^2\ln \frac{R}{l} +2\gamma^2(1-\sqrt{\frac{|y-x|}{l}})& \quad \text{if} \quad |y-x| \leq l, \\
\gamma^2\ln^{+}\frac{R}{|y-x|} & \quad \text{if} \quad |y-x| > l.
\end{cases}
\end{equation*}
One can note the following scale invariance property for $(X_{l}(x))_{x \in \R^2}$: if $\lambda \in ]0,1]$ and $l \leq R$, the two fields $(X_{ \lambda l}(\lambda x))_{|x| \leq R}$ and $(\Omega_{\lambda}+X_{l}(x))_{|x| \leq R}$ have the same law, where $\Omega_{\lambda}$ is a centered Gaussian random variable independent from $(X_{l}(x))_{x \in \R^2}$ and of variance $\gamma^2\ln\frac{1}{\lambda}$.
By taking the limit as $l$ goes to $0$, we get the following scale invariance for $M$: if $\lambda \in ]0,1]$, we have the following identity in law:
\begin{equation}\label{eq:invariance}
(M(\lambda A))_{A \subset B(0,R)}\overset{(Law)}{=}\lambda^2e^{\Omega_{\lambda}-\frac{\gamma^2}{2}\ln\frac{1}{\lambda}}(M(A))_{A \subset B(0,R)}.
\end{equation}
Taking the expectation in (\ref{eq:invariance}) to the power $q \in [0,1]$, we get:
\begin{equation*}
\E[M(B(0,\lambda))^q]=(\frac{\lambda}{R})^{\zeta(q)}\E[M(B(0,R))^q]
\end{equation*}
with:
\begin{equation*}
\zeta(q)=(2+\gamma^2)q-\frac{\gamma^2}{2}q^2.
\end{equation*}

Finally, it is possible to extend naturally the notion of Hausdorff content (and Hausdorff dimension) on a metric space $(X,d)$ to a measurable space $X$ equipped with a measure $\mu$ by :
\begin{equation*}
C^s_H(K)=\inf\left\{\sum_i \mu(B(x_i,r_i))^s;\text{there is a cover of K by balls} \; B(x_i,r_i) \text{with radii }  r_i>0\right\}.
\end{equation*}

With these extensions, we can state the following $2$-dimensional analog to theorem \ref{main}:

\begin{theorem}\label{th:2dim}
Assume that $\gamma^2<4$. Let $K\subset B(0,R)$ be some deterministic and measurable nonempty set and $\delta_0$ its Hausdorff dimension with respect to the Euclidian metric. Then the Hausdorff dimension ${\rm dim}_H^M(K)$ of $K$ with respect to the random measure $M$ coincides $\P$ a.s. with the unique solution $\delta$ in $[0,1] $ of the equation $\delta_0=\zeta(\delta)$.
\end{theorem}

\proof

Just note that, in this setting, the Frostman lemma is unchanged if we define the capacity of $M$ by the following formula:
\begin{equation*}
{\rm Cap}_s(K)=\inf\left\{\Big(\int_{K\times K}(M(x,|y-x|))^{-s}\gamma(dx)\gamma(dy)\big)^{-1};\,\gamma \text{ is a Borel measure such that }\gamma(K)=1\right\}
\end{equation*}

The proof is then a straightforward adaptation of the proof of theorem \ref{main}.

\qed

\subsection{The exponential of the Gaussian Free Field}
In this subsection, as an application of the previous subsection, we prove the KPZ formula for the exponential of the Gaussian Free Field (GFF) in $B(0,R)$: this corresponds in $B(0,R)$ to the gravity measure considered on a 2 dimensional surface in \cite{cf:Da}. The GFF is an important object in Conformal Field theory since it has the conformal invariance property and a spatial Markovian property (see \cite{cf:Sh}). Formally, the GFF (or Euclidian bosonic massless free field) in $B(0,R)$ is a "Gaussian Field" $X$ with covariance given by:
\begin{equation*}
\E[X_F(x)X_F(y)]=G_{R}(x,y),
\end{equation*}
where $G_R$ is the Green function of $B(0,R)$ (see for example chapter 2.4 in \cite{cf:Law} for the definition and main properties). Let the process $B_t$ be Brownian motion starting from
 $x$ under the measure $P^x$ and consider the stopping time $T_R=\text{inf} \{t \geq 0, \; |B_t|=R \}$.
  If we denote $p_R(t,x,y)=P^x(B_{t} \in dy, \; T_R > t)$, we have:
\begin{equation*}
G_{R}(x,y)= \pi \int_{0}^{\infty}p_R(t,x,y)dt.
 \end{equation*}
Note that for each $t>0$, $p_R(t,x,y)$ is a continuous positive and positive definite kernel on $B(0,R)$.
 Therefore, we can define the GFF measure $M_F$ as multiplicative chaos (\cite{cf:Kah}) associated to the
  kernel $\gamma^2G_R$ where $\gamma^2<4$. In this framework, $M_F$ is the almost sure limit (in the space
   of Radon measures) as $l$ goes to $0$ of the measure:
\begin{equation*}
M_{l,F}=e^{X_{l,F}(x)-\frac{1}{2}\E[X_{l,F}(x)^2]}dx
\end{equation*}
where $X_{l,F}$ is a Gaussian field with the following covariance:
\begin{equation*}
\E[X_{l,F}(x)X_{l,F}(y)]=\gamma^2 \pi \int_{l^2}^{+ \infty}p_R(t,x,y)dt.
\end{equation*}

We know have the following analog of theorem \ref{th:2dim}:

\begin{theorem}\label{th:GFF}
Assume that $\gamma^2<4$ and $r<R$. Let $K\subset B(0,r)$ be some deterministic and measurable nonempty set and $\delta_0$ its Hausdorff dimension with respect to the Euclidian metric. Then the Hausdorff dimension ${\rm dim}_H^M(K)$ of $K$ with respect to the random measure $M$ coincides $\P$ a.s. with the unique solution $\delta$ in $[0,1] $ of the equation $\delta_0=\zeta(\delta)$.
\end{theorem}



\section{Proof of Theorem \ref{main}}


\begin{lemma}
Let $x<y\in\R$. If $q\in[0,1]$ then
$$\E[\rho(x,y)^q] \leq C(T,q)|x-y|^{\zeta(q)},$$ where $C(T,q)$ is a positive constant only depending on $T,q$. As a consequence, if $K,\delta,\delta_0$ are defined as in Theorem \ref{main}, then a.s.
$\zeta({\rm dim}_H^\rho(K))\leq \delta_0$.
\end{lemma}

\noindent {\bf Proof.} By stationarity of the measure $M$ and Proposition \ref{propMRM}, we have
\begin{align*}
\E[\rho(x,y)^q]=\E[M([x,y])^q]=\E[M([0,y-x])^q]  =|y-x|^{\zeta(q)}T^{-\zeta(q)}\E[M([0,T])^q].
\end{align*}
So we can choose $C(T,q)= T^{-\zeta(q)}\E[M([0,T])^q]<+\infty$.

Let $\alpha>0$ and $q\in [0,1]$ such that $\zeta(q)>\delta_0$. There exists a covering of $K$ by a countable family $([x_n,y_n])_n$ such that $\sum_n|x_n-y_n|^{\zeta(q)}<\alpha $. Hence
\begin{align*}
\E\Big[\sum_n\rho(x_n,y_n)^q\Big]&=\sum_n\E\big[\rho(x_n,y_n)^q\big]\leq C(T,q)
\sum_n|y_n-x_n|^{\zeta(q)}\leq C(T,q)\alpha.
\end{align*}
By the Markov inequality, $\P\big(\sum_n\rho(x_n,y_n)^q\leq C(T,q)\sqrt{\alpha}\big)\geq 1-\sqrt{\alpha}$.  Put in other words,with probability $1-\sqrt{\alpha} $, we have a covering of K with balls whose $\rho$-radii satisfy $\sum_n\rho(x_n,y_n)^q\leq  C(T,q)\sqrt{\alpha}$. Thus $q\geq{\rm dim}_H^\rho(K) $ a.s. and the lemma follows.\qed

\begin{proposition}\label{prop_frost}
Let $K,\delta,\delta_0,{\rm dim}_H^\rho(K)$ be as in Theorem \ref{main} and let $q\in[0,1]$ be such that $\zeta(q)<\delta_0 $. Then a.s. $q\leq {\rm dim}_H^\rho(K) $, that is $\delta_0\leq \zeta({\rm dim}_H^\rho(K)) $.
\end{proposition}

\noindent {\bf Proof.} Since $\zeta(q)<\delta_0 $, by the Frostman Lemma, there is a Borel probability measure $\gamma_0 $ supported by $K$ such that $\gamma_0(K)=1$ and
$$\int_{[0,T]^2}|x-y|^{-\zeta(q)}\,\gamma_0(dx)\,\gamma_0(dy)<+\infty. $$
Let us define, for any $0<l<T$, the measure on $[0,T]$:
$$\nu_l(dr)=e^{q \omega_l(r)-\psi(q)(\ln(T/l)+1)}\,\gamma_0(dr) $$ and its associated metric on $\R$:
$$\forall x,y\in \R,\quad \rho_l(x,y)=\nu_l([x,y]). $$
We now investigate the quantity:
\begin{align*}
\phi(l,\gamma_0)&\equiv\E\Big[\int _{[0,T]^2}\rho_l(x,y)^{-q}\,\nu_l(dx)\,\nu_l(dy)\Big]\\
&=\int _{[0,T]^2}\E\Big[\rho_l(x,y)^{-q}e^{q\omega_l(x)+q\omega_l(y)-2\psi(q)(\ln(T/l)+1)}\Big]\,\gamma_0(dx)\gamma_0(dy)\\&=2\int _{y\geq x}\E\Big[\rho_l(0,y-x)^{-q}e^{q\omega_l(0)+q\omega_l(y-x)-2\psi(q)(\ln(T/l)+1)}\Big]\,\gamma_0(dx)\gamma_0(dy)
\end{align*}
by stationarity of the process  $\omega_l$. To this purpose, we split the above integral in two terms as
\begin{align*}
\phi(l,\gamma_0)=& 2\int _{0 \leq y-x< l}\E\Big[\rho_l(0,y-x)^{-q}e^{q\omega_l(0)+q\omega_l(y-x)-2\psi(q)(\ln(T/l)+1)}\Big]\,\gamma_0(dx)\gamma_0(dy)\\&+2\int _{ y-x\geq l}\E\Big[\rho_l(0,y-x)^{-q}e^{q\omega_l(0)+q\omega_l(y-x)-2\psi(q)(\ln(T/l)+1)}\Big]\,\gamma_0(dx)\gamma_0(dy)\\
\equiv & \phi_1(l,\gamma_0)+\phi_2(l,\gamma_0).
\end{align*}

We first estimate $ \phi_1(l,\gamma_0)$. Using the Jensen inequality and the decrease of the mapping $x\mapsto x^{-q}$ yields
\begin{align*}
&\phi_1(l,\gamma_0)\\
=&2\int _{0 \leq y-x< l}\E\Big[\Big(\int_0^{y-x}e^{\omega_l(r)\,dr}\Big)^{-q}e^{q\omega_l(0)+q\omega_l(y-x)-2\psi(q)(\ln(T/l)+1)}\Big]\,\gamma_0(dx)\gamma_0(dy)\\
=& \frac{2 e^{-2\psi(q)}l^{2\psi(q)}}{T^{2\psi(q)}}\int _{0 \leq y-x< l}\E\Big[\Big(\int_0^{y-x}e^{\omega_l(r)-\omega_l(0)-\omega_l(y-x)}\,dr\Big)^{-q}\Big]\,\gamma_0(dx)\gamma_0(dy)\\
\leq & \int _{0 \leq y-x< l}\frac{2 e^{-2\psi(q)}l^{2\psi(q)}}{T^{2\psi(q)}|y-x|^{q}}\E\Big[e^{\int_0^{y-x}(q\omega_l(0)+q\omega_l(y-x)-q\omega_l(r))\frac{dr}{y-x}}\Big]\,\gamma_0(dx)\gamma_0(dy).
\end{align*}
Given $0\leq x<y\leq T$ such that $y-x<l$, define $A_l^i\equiv A_l(0)\cap A_l(y-x) \not =\emptyset$. Each cone-like subset $A_l(r)$ ($0\leq r\leq y-x$) can be split into three terms as $A_l(r)=A^{g}_l(r)\cup A_l^i\cup A_l^d(r)$, where $A^{g}_l(r)$ (resp. $A^{d}_l(r)$) denotes the part of $A_l(r)$ located on the left (resp. right) of $A^i_l$. It is worth emphasizing that:
 $$
 (\omega^d_l(r))_{0\leq r \leq y-x}=(\mu(A_l^d(y-x) \setminus A_l^d(y-x-r))-\psi'(0)\theta(A_l^d(y-x) \setminus A_l^d(y-x-r)))_{0\leq r \leq y-x}
 $$
 is a right-continuous martingale, as well as $(\omega^g_l(r))_{0\leq r \leq y-x}$ where:
 $$
 \omega^g_l(r)=\mu(A^{g}_l(0) \setminus A^{g}_l(r))-\psi'(0)\theta(A^{g}_l(0) \setminus A^{g}_l(r)).
 $$
 By using the fact that $\psi'(0)<0$, we get:
\begin{align*}
q\omega_l(0)+q\omega_l(y-x)-q\omega_l(r) & = q\omega^i_l+q\mu(A_l^d(y-x) \setminus A_l^d(r))+q\mu(A^{g}_l(0) \setminus A^{g}_l(r))   \\
& \leq q\omega^i_l+q\omega_l^d(y-x-r)+q \omega^g_l(r). \\
\end{align*}
 Since $(\omega^d_l(r))_{r}$, $(\omega^g_l(r))_{r}$ and $w_l^i=\mu(A_l^i)$ are independent, the last expression is estimated as:
\begin{align*}
&\phi_1(l,\gamma_0)\\
\leq& \int _{0 \leq y-x< l}\frac{2 e^{-2\psi(q)}l^{2\psi(q)}}{T^{2\psi(q)}|y-x|^{q}}\E[e^{q\omega^i_l}]\E[\sup_{0\leq r \leq y-x}e^{q\omega_l^d(y-x-r)}]\E[\sup_{0\leq r \leq y-x}e^{q\omega^g_l(r)}]\,\gamma_0(dx)\gamma_0(dy)\\
\leq& \int _{0 \leq y-x< l}\frac{2 C_q^2 e^{-2\psi(q)}l^{2\psi(q)}}{T^{2\psi(q)}|y-x|^{q}}\E[e^{q\omega^i_l}]\E[e^{q\omega_l^d(y-x)}]\E[e^{q\omega^g_l(y-x)}]\,\gamma_0(dx)\gamma_0(dy),
\end{align*}
the last inequality resulting from the Doob inequality applied to the function $x \rightarrow e^x$ ($C_q$ is a constant only depending on $q$).
It remains to compute $ \theta(A_l^i)$, $\theta(A_l^g(0)) $ and $\theta(A_l^d(y-x))$. It is plain to see that
$$  \theta(A_l^i)=\ln(T/l)+1-(y-x)/l,\quad \theta(A_l^d(y-x))=\theta(A_l^g(0)) =(y-x)/l,$$ in such a way that
(we use that $\psi(q)<0$ for all $q$ in $]0,1[$):
\begin{align}\label{ineg1}
\phi_1(l,\gamma_0)\leq  & \int _{0 \leq y-x< l}\frac{2 C_q^2 e^{-2\psi(q)}l^{2\psi(q)}}{T^{2\psi(q)}|y-x|^{q}}e^{\psi(q)\big(\ln(T/l)+1+(y-x)/l\big)}e^{2(\psi(q)\frac{y-x}{l}-\psi'(0)\frac{y-x}{l})}\gamma_0(dx)\gamma_0(dy)\\
\leq  & 2e^{-2\psi'(0)}C_q^2(eT)^{-\psi(q)}\int _{0 \leq y-x< l}   \frac{1}{|y-x|^{\zeta(q)}}    \gamma_0(dx)\gamma_0(dy).\nonumber
\end{align}

Let us now focus on $ \phi_2(l,\gamma_0)$. In what follows, we make a change of variable $u=Tr/(y-x)$:
\begin{align*}
 & \phi_2(l,\gamma_0)\\
&=2\int _{y-x\geq l}\E\Big[\frac{e^{q\omega_l(0)+q\omega_l(y-x)-2\psi(q)(\ln(T/l)+1)}}{\Big(\int_0^{y-x}e^{\omega_l(r)\,dr}\Big)^{q}}\Big]\,\gamma_0(dx)\gamma_0(dy)\\
&=  \int _{y-x\geq l}\frac{2T^q}{|y-x|^q}\E\Big[\frac{e^{q\omega_l(0)+q\omega_l(y-x)-2\psi(q)(\ln(T/l)+1)}}{\Big(\int_0^Te^{\omega_l((y-x)uT^{-1})}\,du\Big)^{q}}\Big]\,\gamma_0(dx)\gamma_0(dy)
\end{align*}
We remind the reader of the following property: the process $(\omega_{l' \alpha}(\alpha t))_{0\leq t \leq T} $ has the same law as the process $(\Omega_{\alpha}+\omega_{l' }(t))_{0\leq t \leq T} $, where $\alpha\in]0,1]$, $l' \leq T$ and $\Omega_{\alpha}$ is an infinitely divisible random variable independent from the process $(\omega_{l'}(t))_{0\leq t \leq T} $ such that $ \E[e^{iq \Omega_{\alpha}}]=\alpha^{-\varphi(q)}$. In particular, choosing $l'=lT/(y-x) $ and $\alpha=(y-x)/T $, the process $\big(\omega_{l}\big((y-x) t/T\big)\big)_{0\leq t \leq T} $ has the same law as the process $(\Omega_{(y-x)/T}+\omega_{lT/(y-x) }(t))_{0\leq t \leq T} $. Plugging this relation into the above estimate of $\phi_2(l,\gamma_0)$ yields
\begin{align*}
 & \phi_2(l,\gamma_0)  \\
 &=  \int _{y-x\geq l}\frac{2T^q}{|y-x|^q}\E\Big[\frac{e^{q\Omega_{(y-x)/T}+q\omega_{\frac{lT}{y-x}}(0)+q\omega_{\frac{lT}{y-x}}(T)-2\psi(q)(\ln(T/l)+1)}}{\Big(\int_0^Te^{\omega_{\frac{lT}{y-x}}(u)}\,du\Big)^{q}}\Big]\,\gamma_0(dx)\gamma_0(dy)  \\
&=  \int _{y-x\geq l}\frac{2T^{\zeta(q)}}{|y-x|^{\zeta(q)}}\E\Big[\frac{e^{q\omega_{\frac{lT}{y-x}}(0)+q\omega_{\frac{lT}{y-x}}(T)-2\psi(q)(\ln(\frac{y-x}{l})+1)}}{\Big(\int_0^Te^{\omega_{\frac{lT}{y-x}}(u)}\,du\Big)^{q}}\Big]\,\gamma_0(dx)\gamma_0(dy)  \\
\end{align*}
 Thus it just remains to show that there exists $C >0$ such that for all $l'$ in $[0,T]$:
 \begin{equation*}
 \E\Big[\frac{e^{q\omega_{l'}(0)+q\omega_{l'}(T)-2\psi(q)(\ln(T/l')+1)}}{\Big(\int_0^Te^{\omega_{l'}(u)}\,du\Big)^{q}}\Big] \leq C
 \end{equation*}
In the above inequality, we will restrict to the (non obvious) case $l' \in [0,T/4]$. We have:
\begin{align*}
& \E\Big[\frac{e^{q\omega_{l'}(0)+q\omega_{l'}(T)-2\psi(q)(\ln(T/l')+1)}}{\Big(\int_0^Te^{\omega_{l'}(u)}\,du\Big)^{q}}\Big]  \\
& \leq \E\Big[\frac{e^{q\omega_{l'}(0)+q\omega_{l'}(T)-2\psi(q)(\ln(T/l')+1)}}{\Big(\int_{T/4}^{3T/4}e^{\omega_{l'}(u)}\,du\Big)^{q}}\Big]
\end{align*}
It is worth mentioning that the sets $A_{l'}(0)$, $A_{l'}(T)$ are disjoint. We then define
\begin{align*}
B^g_{l'} & =A_{l'}(0)\setminus A_{l'}(T/4)  \\
B^d_{l'} & =A_{l'}(T)\setminus A_{l'}(3T/4)  \\
\end{align*}
We stress that for any $u$ in $[T/4,3T/4]$:
$$
A_{l'}(u)\cap B^g_{l'}=\emptyset,\quad
A_{l'}(u)\cap  B^d_{l'}=\emptyset
$$
Using the relation $\theta(B^g_{l'})=\theta(B^g_{l'})=\ln(T/l')+1-\ln(4)$ and the independence of
$\mu(B^g_{l'})$, $\mu(B^d_{l'})$, $(\mu(A_{l'}(u)))_{T/4 \leq u \leq 3T/4}$, we get:
\begin{align*}
 & \E\Big[\frac{e^{q\omega_{l'}(0)+q\omega_{l'}(T)-2\psi(q)(\ln(T/l')+1)}}{\Big(\int_{T/4}^{3T/4}e^{\omega_{l'}(u)}\,du\Big)^{q}}\Big]  \\
 & = e^{-2\psi(q)(\ln(T/l')+1)} \E\Big[ e^{q\mu(B^g_{l'})}\Big] \E\Big[  e^{q\mu(B^d_{l'})}\Big]  \E\Big[\frac{e^{q\mu(A_{l'}(0) \cap  A_{l'}(T/4))+q\mu(A_{l'}(T) \cap  A_{l'}(3T/4))}}{\Big(\int_{T/4}^{3T/4}e^{\omega_{l'}(u)}\,du\Big)^{q}}\Big]  \\
& = e^{-2\ln(4)\psi(q)}\E\Big[\frac{e^{q\mu(A_{l'}(0) \cap  A_{l'}(T/4))+q\mu(A_{l'}(T) \cap  A_{l'}(3T/4))}}{\Big(\int_{T/4}^{3T/4}e^{\omega_{l'}(u)}\,du\Big)^{q}}\Big]  \\
\end{align*}

Let us denote $\mathcal{A}_{l'}^g(u),\mathcal{A}_{l'}^d(u)$ the following sets for $u \in [T/4,3T/4]$:
\begin{align*}
\mathcal{A}_{l'}^g(u)  & =   (A_{l'}(0) \cap  A_{l'}(u)) \setminus  A_{l'}(3T/4) \\
\mathcal{A}_{l'}^d(u)  & =    (A_{l'}(T) \cap  A_{l'}(u)) \setminus  A_{l'}(T/4)  \\
\end{align*}

We have the following decompositions:
\begin{align*}
\mu\Big( A_{l'}(0) \cap  A_{l'}(T/4) \Big) & =\mu\Big( \mathcal{A}_{l'}^g(T/4)  \Big)
+ \mu\Big( A_{l'}(0) \cap   A_{l'}(3T/4) \Big), \\
\mu\Big( A_{l'}(T) \cap  A_{l'}(3T/4) \Big) & =\mu\Big(  \mathcal{A}_{l'}^d(3T/4) \Big)
+ \mu\Big( A_{l'}(T)  \cap   A_{l'}(T/4) \Big). \\
\end{align*}
We also have for all $u$ in $[T/4,3T/4]$:
\begin{align*}
\mu \Big( A_{l'}(u) \Big) & = \mu\Big( \mathcal{A}_{l'}^g(u)  \Big)  \\
& + \mu\Big( A_{l'}(0) \cap   A_{l'}(3T/4) \Big)  \\
& +\mu\Big( \mathcal{A}_{l'}^d(u)  \Big)  \\
& + \mu\Big( A_{l'}(T)  \cap   A_{l'}(T/4) \Big)  \\
& +  \mu\Big( A_{l'}(u)\setminus (A_{l'}(0)\cup A_{l'}(T) \Big).  \\
\end{align*}

Therefore, we get:

\begin{align*}
 & \E\Big[\frac{e^{q\omega_{l'}(0)+q\omega_{l'}(T)-2\psi(q)(\ln(T/l')+1)}}{\Big(\int_{T/4}^{3T/4}e^{\omega_{l'}(u)}\,du\Big)^{q}}\Big]  \\
 & = e^{-2\ln(4)\psi(q)}\E\Big[\frac{e^{q\mu( A_{l'}(0) \cap  A_{l'}(T/4) )+q\mu( A_{l'}(T) \cap  A_{l'}(3T/4))}}{\Big(\int_{T/4}^{3T/4}e^{\omega_{l'}(u)}\,du\Big)^{q}}\Big]  \\
& = e^{-2\ln(4)\psi(q)}\E\Big[\frac{e^{q\mu(  \mathcal{A}_{l'}^g(T/4))+q\mu(  \mathcal{A}_{l'}^d(3T/4))}}{\Big(\int_{T/4}^{3T/4}e^{ \mu(  \mathcal{A}_{l'}^g(u) )
+\mu(  \mathcal{A}_{l'}^d(u)   )
+  \mu( A_{l'}(u)\setminus (A_{l'}(0)\cup A_{l'}(T)) )}\,du\Big)^{q}}\Big]  \\
& \leq e^{-2\ln(4)\psi(q)} \E \Big[ e^{q\mu( \mathcal{A}_{l'}^g(T/4) )- q\inf_u\mu(  \mathcal{A}_{l'}^g(u) )}\Big] \\
& \times  \E \Big[ e^{q\mu(  \mathcal{A}_{l'}^d(3T/4)   )- q\inf_u\mu(   \mathcal{A}_{l'}^d(u)  )}\Big]  \E \Big[    \frac{1}{\Big(\int_{T/4}^{3T/4}e^{\mu( A_{l'}(u)\setminus (A_{l'}(0)\cup A_{l'}(T) )}\,du\Big)^{q}}   \Big]  \\
& = e^{-2\ln(4)\psi(q)} \E \Big[ e^{q\sup_{u}(\mu( \mathcal{A}_{l'}^g(T/4)  )- \mu(  \mathcal{A}_{l'}^g(u) ))}\Big] \\
& \times  \E \Big[ e^{q\sup_u( \mu(  \mathcal{A}_{l'}^d(3T/4)   )- \mu(   \mathcal{A}_{l'}^d(u) ))}\Big]  \E \Big[    \frac{1}{\Big(\int_{T/4}^{3T/4}e^{\mu( A_{l'}(u)\setminus (A_{l'}(0)\cup A_{l'}(T) )}\,du\Big)^{q}}   \Big]  \\
& = e^{-2\ln(4)\psi(q)} \E \Big[ e^{q\sup_{u}(\mu( \mathcal{A}_{l'}^g(T/4) \setminus \mathcal{A}_{l'}^g(u) ))}\Big] \\
& \times  \E \Big[ e^{q\sup_u( \mu(  \mathcal{A}_{l'}^d(3T/4)  \setminus  \mathcal{A}_{l'}^d(u) ))}\Big]  \E \Big[    \frac{1}{\Big(\int_{T/4}^{3T/4}e^{\mu( A_{l'}(u)\setminus (A_{l'}(0)\cup A_{l'}(T) )}\,du\Big)^{q}}   \Big]  \\
\end{align*}

The process
\begin{equation*}
\mu( \mathcal{A}_{l'}^g(T/4) \setminus \mathcal{A}_{l'}^g(u) )-\psi'(0)\theta(\mathcal{A}_{l'}^g(T/4) \setminus \mathcal{A}_{l'}^g(u))
\end{equation*}
 is a martingale for $u$ in $[T/4,3T,4]$ and we have $\theta(\mathcal{A}_{l'}^g(T/4))$ bounded independently from $l'$. By applying Doob's inequality, there exists some constant $C>0$ independent from $l'$ such that:
 \begin{equation*}
  \E \Big[ e^{q\sup_{u}(\mu( \mathcal{A}_{l'}^g(T/4) \setminus \mathcal{A}_{l'}^g(u) ))}\Big]  \leq C.
 \end{equation*}
Similarly, we have:
 \begin{equation*}
 \E \Big[ e^{q\sup_u( \mu(  \mathcal{A}_{l'}^d(3T/4)  \setminus  \mathcal{A}_{l'}^d(u) ))}\Big] \leq C
 \end{equation*}
Therefore, we get:

\begin{align*}
 & \E\Big[\frac{e^{q\omega_{l'}(0)+q\omega_{l'}(T)-2\psi(q)(\ln(T/l')+1)}}{\Big(\int_{T/4}^{3T/4}e^{\omega_{l'}(u)}\,du\Big)^{q}}\Big]  \\
& \leq C \E \Big[    \frac{1}{\Big(\int_{T/4}^{3T/4}e^{\mu( A_{l'}(u)\setminus (A_{l'}(0)\cup A_{l'}(T) )}\,du\Big)^{q}}   \Big]  \\
\end{align*}
Since $\psi(-q)<\infty$, by using the same argument than the proof of theorem 3 (Moments of negative orders) in \cite{Bar}, one can show that:
\begin{equation*}
\sup_{l'} \E \Big[    \frac{1}{\Big(\int_{T/4}^{3T/4}e^{\mu\Big( A_{l'}(u)\setminus (A_{l'}(0)\cup A_{l'}(T) \Big)}\,du\Big)^{q}}   \Big] < \infty.
\end{equation*}

To sum up, gathering the estimates of $\phi_1(l,\gamma_0)$ and $\phi_2(l,\gamma_0)$, we have proved the existence of some constant $C>0$ such that:
$$ \phi(l,\gamma_0)\leq C \int_{[0,T]^2}\frac{1}{|y-x|^{\zeta(q)}}\gamma_0(dx)\gamma_0(dy)<+\infty.$$
Let us now define the measure $\nu(dt)=\lim_{l\to 0^+}\nu_l(dt)$ (see Lemma \ref{weaknu} below). From Lemma \ref{weaknu} and the Fatou lemma, we obtain
\begin{align*}
\E\Big[\int _{[0,T]^2}\rho(x,y)^{-q}\,\nu(dx)\,\nu(dy)\Big] & \leq \E\Big[\liminf_{l\to 0^+}\int _{[0,T]^2}\rho_l(x,y)^{-q}\,\nu_l(dx)\,\nu_l(dy)\Big]\\
& \leq \liminf_{l\to 0^+}\E\Big[\int _{[0,T]^2}\rho_l(x,y)^{-q}\,\nu_l(dx)\,\nu_l(dy)\Big]\\
&\leq C\int_{[0,T]^2}\frac{1}{|y-x|^{\zeta(q)}}\gamma_0(dx)\gamma_0(dy)<+\infty.
\end{align*}
As a consequence, $\P$ a.s. the integral $\int _{[0,T]^2}\rho(x,y)^{-q}\,\nu(dx)\,\nu(dy)$ is finite. We complete the proof with the Frostman Lemma.\qed

\begin{lemma}\label{weaknu}
Assume that we are given $q\in [0,1]$ such that
$$\int_{[0,T]^2}\frac{\gamma_0(dx)\gamma_0(dy)}{|y-x|^{\zeta(q)}} <+\infty.$$
We consider, for any $l>0$, the measure on $[0,T]$:
$$\nu_l(dt)=e^{q\omega_l(t)-\psi(q)\big(\ln(T/l)+1\big)}\,\gamma_0(dt) .$$
Then the weak limit (in the sense of measures) $$\nu(dt)=\lim_{l\to 0^+}\nu_l(dt) $$ exists $\P$-a.s., is finite, supported by $K$ $\P$-a.s., and we have
$$\int _{[0,T]^2}\rho(x,y)^{-q}\,\nu(dx)\,\nu(dy)\leq \liminf_{l\to 0^+}\int _{[0,T]^2}\rho_l(x,y)^{-q}\,\nu_l(dx)\,\nu_l(dy).$$
\end{lemma}

\noindent {\bf Proof.} According to the proof of Proposition $\ref{prop_frost} $, we have
$$\phi(l,\gamma_0)\leq C \int_{[0,T]^2}\frac{\gamma_0(dx)\gamma_0(dy)}{|y-x|^{\zeta(q)}} <+\infty.$$ Furthermore, $ \rho_l(x,y)\leq \rho_l(0,T)$ for any $0\leq x\leq y \leq T$, in such a way that
$$\E[\nu_l(A)^2\rho_l(0,T)^{-\zeta(q)}]\leq \phi(l,\gamma_0)\leq C\int_{[0,T]^2}\frac{\gamma_0(dx)\gamma_0(dy)}{|y-x|^{\zeta(q)}} <+\infty$$
for any Lebesgue measurable subset $A$ of $[0,T]$. Moreover, if the Lebesgue measure of $A$ is strictly positive then the H\"older inequality yields
\begin{equation}\label{holder}
\begin{split}
\E[\nu_l(A)^{2/(1+\zeta(q))}]\leq & \E[\nu_l(A)^{2}M_l([0,T])^{-\zeta(q)}]^{1/(1+\zeta(q))}\E[M_l([0,T])]^{\zeta(q)/(1+\zeta(q))}\\
& \leq C'\int_{[0,T]^2}\frac{\gamma_0(dx)\gamma_0(dy)}{|y-x|^{\zeta(q)}} <+\infty.
\end{split}
\end{equation}
We remind the reader that $(\nu_l(A))_l$ is martingale for any Lebesgue measurable subset $A$ of $[0,T]$.  From \eqref{holder}, this martingale is bounded in $L^{1+\epsilon}$ for some $\epsilon>0$. As a consequence, it converges $\P$-a.s. towards a limit denoted by $\nu(A)$ as $l\to 0$. It is readily seen that $\nu$ is a measure on $[0,T]$ $\P$-a.s. Since $\nu_l(K^c)=0$, it is clear that $\nu(K^c)=0$ $\P$-a.s.

Finally, $\E[\nu([0,T])]=\lim_{l\to 0}\E[\nu_l([0,T])]=\gamma_0([0,T])\geq 1$. Moreover $\{\nu([0,T])>0\}$ is an event of the asymptotic $\sigma$-field generated by the random variables $(\nu_l(A))_l$ and has therefore probability 0 or 1. As a consequence,  the event $\{\nu([0,T])>0\}$ has probability 1.

The last inequality of the lemma results from Lemma \ref{unifconv} below and the weak convergence of measures.\qed

\begin{lemma}\label{unifconv}
$\P$ a.s., the metric $(\rho_l)_l$ uniformly converges towards the metric $\rho$ as $l\to 0$, that is
$$\P \text{ a.s.},\quad \lim_{l\to 0}\sup_{0\leq x\leq y \leq T}|\rho_l(x,y)-\rho(x,y)|=0.$$
\end{lemma}

\noindent {\bf Proof.} The mapping $x\mapsto \rho(0,x)$ is continuous because of the non-degeneracy of $\rho$ (see Proposition \ref{propMRM}). Moreover, for each $l>0$, the mapping $x\mapsto \rho_l(0,x)$ is increasing and the sequence $(\rho_l(0,x)$ converges pointwise $\P$ a.s. towards $\rho(0,x)$ (see Definition \ref{wcmrm}). The uniform convergence then results from the Dini theorem.\qed

\section{Proof of theorem \ref{th:GFF}}

Let $r<R$. We choose $\delta>0$ such that $r+\delta<R$. With the notations of section 5.1, one can see that there exists two positive constants $c_{r,\delta},C_{r,\delta}$ such that for all $x,y \in B(0,r+\delta)$ we have (independently of $l$ and $x,y$):
\begin{equation}\label{eq:encadrement}
\E[X_{l}(x)X_{l}(y)] +c_{r,\delta} \leq \E[X_{l,F}(x)X_{l,F}(y)]  \leq  \E[X_{l}(x)X_{l}(y)]+C_{r,\delta}.
\end{equation}

\emph{Proof of}: $\zeta( \text{dim}_H^{M_{F}} (K))  \leq \text{dim}_H(K)$.

The inequality (\ref{eq:encadrement}) and the classical corollary 6.2 in \cite{cf:RoVa} imply the existence for $q \in [0,1]$ of $C_{q,r,\delta}>0$ such that:
\begin{equation*}
\forall B(x_i,r_i) \subset  B(0,r+\delta), \quad \E[M_{F}(B(x_i,r_i))^q] \leq C_{q,r,\delta}  r_i^{\zeta(q)}.
\end{equation*}
We conclude by using the same argument than in the proof of theorem \ref{main}.

\emph{Proof of}: $\zeta( \text{dim}_H^{M_{F}} (K))  \geq \text{dim}_H(K)$.

Suppose $\zeta(q) < \text{dim}_H(K)$. Following the notations of section 6 (proof of theorem \ref{main}), we consider a measure $\gamma_0 $ supported by $K$ such that $\gamma_0(K)=1$ and
$$\int_{[0,T]^2}|x-y|^{-\zeta(q)}\,\gamma_0(dx)\,\gamma_0(dy)<+\infty. $$
The inequality (\ref{eq:encadrement}) and the classical corollary 6.2 in \cite{cf:RoVa} imply the existence of some constant $C_{q,r,\delta}>0$ such that for all $x,y \in B(0,r)$ with $|y-x|\leq \delta$:
\begin{align*}
& \E\Big[(\int_{B(x,|y-x|)}e^{X_{l,F}(z)-\frac{1}{2}\E[X_{l,F}(z)^2]}dz)^{-q}e^{qX_{l,F}(x)+qX_{l,F}(y)-\frac{q^2}{2}\E[X_{l,F}(x)^2]-\frac{q^2}{2}\E[X_{l,F}(y)^2]}\Big]  \\
& \leq C_{q,r,\delta}  \E\Big[(\int_{B(x,|y-x|)}e^{X_{l}(z)-\frac{1}{2}\E[X_{l}(z)^2]}dz)^{-q}e^{qX_{l}(x)+qX_{l}(y)-\frac{q^2}{2}\E[X_{l}(x)^2]-\frac{q^2}{2}\E[X_{l}(y)^2]}\Big]  \\
\end{align*}

Taking the limit as $l$ goes to $0$, this implies:
\begin{align*}
& \; \underset{l \to \infty }{\underline{\lim}} \int_{|y-x| \leq \delta} \gamma_0(dx)\gamma_0(dy)\E\Big[(\int_{B(x,|y-x|)}e^{X_{l,F}(z)-\frac{1}{2}\E[X_{l,F}(z)^2]}dz)^{-q}e^{qX_{l,F}(x)+qX_{l,F}(y)-\frac{q^2}{2}\E[X_{l,F}(x)^2]-\frac{q^2}{2}\E[X_{l,F}(y)^2]}\Big]  \\
& \leq  C_{q,r,\delta} \underset{l \to \infty }{\underline{\lim}} \int_{|y-x| \leq \delta} \gamma_0(dx)\gamma_0(dy) \E\Big[(\int_{B(x,|y-x|)}e^{X_{l}(z)-\frac{1}{2}\E[X_{l}(z)^2]}dz)^{-q}e^{qX_{l}(x)+qX_{l}(y)-\frac{q^2}{2}\E[X_{l}(x)^2]-\frac{q^2}{2}\E[X_{l}(y)^2]}\Big]  \\
& < \infty.
\end{align*}

We remind that the second inequality above results from a straightforward adaptation to the $2$ dimensional case  of the proof of theorem \ref{main} (in the log normal case).

We then conclude by using the same argument than in the proof of theorem \ref{main}.


\end{document}